\documentstyle[12pt]{article}

\begin{document}

\title{Amenability of Hopf $C^*$-algebras}
\author{Chi-Keung Ng}
\date{\today}
\maketitle

\begin{abstract}
In this paper, we will define and study amenability of Hopf $C^*$-algebras
as well as that of Fourier duality.
\end{abstract}

\noindent{\small 1991 AMS Mathematics Classification number: 46L55, 46L05}

\medskip

\medskip

\medskip

\noindent {\large 0. Introduction}

\medskip

\medskip

Amenability is an interesting and well studied subject in the theory of
locally compact groups. There are a lot of nice properties associated with
it. Since Hopf $C^*$-algebras are generalisations of locally compact groups,
it is natural to ask whether we can define amenability for Hopf $C^*$%
-algebras and generalise some results in this situation. There are already
some generalisations of amenability to locally compact semi-groups (see e.g.
[17] or [14]), Banach algebras (see e.g. [15]) as well as operator algebras
(see e.g. [16]). Even though these are interesting and well developed
subjects by themselves, they do not seem to be appropriate approaches for
general Hopf $C^*$-algebras.

\medskip

The aims of this paper are to find good definitions for amenability of Hopf $%
C^*$-algebras and study their properties. Our approach is based on the
properties of (generalised) Fourier algebras and part of our consideration
is on the same line as that of Kac algebras as in [5].

\medskip

It is well known that a locally compact group $G$ is amenable if and only if 
$C^*_r(G) = C^*(G)$ which is the case if and only if the Fourier algebra $%
A(G)$ separates points of $C^*(G)$ (see also [8, 4.1]). It is also well
known that $G$ is amenable if and only if $A(G)$ has a bounded approximate
identity (see [8,4.1]). Moreover, it was proved in [9, Theorem 1] that the
amenability of $G$ is equivalent to the fact that $B(G)=M(A(G))$ (where $B(G)
$ is the Fourier-Stieltjes algebra and $M(A(G))$ is the space $M(A(G),A(G))$
defined in Definition 1.11).

\medskip

Note that the above characterisations all centre around the Fourier algebra $%
A(G)$. In [12, \S 5], we defined the Fourier algebra $A_S$ of a general Hopf 
$C^*$-algebra $S$. We showed that it is a nice generalisation of $A(G)$ and
in most of the interesting cases, it is easy to obtain (see [12, 5.9]). It
also helps us to formulate and prove a uniqueness type theorem of Fourier
duality.

\medskip

Therefore, it is natural to define coamenability of $S$ in one of the
following ways:

\noindent (a) $A_S$ has a bounded approximate identity;

\noindent (b) $A_S$ separates points of $S$ and $S$ is counital;

\noindent (c) $S^* \cong M(A_S)$.

\noindent Note that if $S=C_0(G)$ (which should be regarded as an 
coamenable Hopf $C^*$-algebra), these three conditions are satisfied 
automatically. We
particularly consider these conditions since they can be easily and directly
translated into the case of Hopf $C^*$-algebras without considering the dual
Hopf $C^*$-algebras.

\medskip

It is natural to ask whether these conditions are the same or at least
related to one another. We will give an answer to this question in section
2. Moreover, we will also study some properties and consequences of them.

\medskip

Another natural question is what may happen if we consider amenability in
the duality framework. More precisely, if $(T,X,S)$ is a Fourier duality and 
$T$ is coamenable in any one of the above senses, what property $S$ should
have. Because of the uniqueness of Kac-Fourier duality, it is believed that
there should be a way to characterise coamenability of $T$ in terms of
property of $S$. It turns out to be the amenability of $S$ (as defined in
Definition 3.1(b)).

\medskip

Note that in the case of a locally compact group $G$, we can characterise
the amenability of $G$ by the existence of an invariant mean on $L^\infty(G)$%
. We do not know whether in the Hopf $C^*$-algebra case such a definition
will be equivalent to the previous ones. Nevertheless, we obtain some
relations concerning invariant means in Propositions 3.6 and 3.12.

\medskip

Now, we give an outline of the paper.

\medskip

In section 1, we will give some basic definitions and simple results on Hopf 
$C^*$-algebras and their dualities. We will also recall some basic notions
and facts on Banach algebras which are needed in the following sections.

\medskip

In section 2, we will give three definitions of coamenability of Hopf $C^*$%
-algebras according to the idea above. We will show that the first
definition is the strongest amongst the three and the second one is stronger
than the third one in the 2-sided case. We will also give some interesting
implications of the above coamenability conditions.

\medskip

In section 3, we will study amenability in duality framework (i.e. the
Fourier duality). We will define and study amenable covariant
representations as well as amenable Hopf $C^*$-algebras. In particular, we
will give a list of equivalent conditions for the amenability of covariant
representations as in [5, 2.4]. We will also study nuclearity of crossed
products as in [11]. As a consequence, if $S$ is amenable, then the dual
object $\hat S_p$ is a nuclear $C^*$-algebra. We will apply these results to
Kac-Fourier dualities and show that the three coamenabilities in section 2 are
essentially the same in this case and coincide with the amenability of the
(reduced) dual Hopf $C^*$-algebra.

\medskip

\medskip

\medskip

\noindent {\large 1. Preliminary results}

\medskip

\medskip

We begin by recalling the definition of Hopf $C^*$-algebras. Note that all
Hopf $C^*$-algebras considered in this paper are ``saturated'' in the
following sense.

\medskip

\noindent {\bf Definition 1.1:} Let $A$ and $S$ be $C^*$-algebras.

\noindent (a) $S$ is said to be a {\it Hopf $C^*$-algebra} if there is a
non-degenerate $*$-homomorphism $\delta_S$ from $S$ to $M(S\otimes S)$ such
that\vspace{-3mm}:

\begin{enumerate}
\item[i.]  $(\delta_S\otimes {\rm id})\delta_S = ({\rm id}\otimes
\delta_S)\delta_S$\vspace{-3mm};

\item[ii.]  $\delta_S(1\otimes S), \delta_S(S\otimes 1)\subseteq S\otimes S$%
\vspace{-3mm}.
\end{enumerate}

Moreover, $S$ is said to be {\it saturated} if the linear spans of the sets $%
\{\delta_S(s)(1\otimes t): s,t\in S\}$ and $\{\delta_S(s)(t\otimes 1):
s,t\in S\}$ are dense in $S\otimes S$.

\medskip \noindent (b) A non-degenerate $*$-homomorphism $\epsilon$ from $A$
to $M(A\otimes S)$ is said to be a {\it coaction} if $(\epsilon\otimes {\rm %
id})\epsilon = ({\rm id}\otimes \delta_S)\epsilon$ and $N=\{\epsilon(a)(1%
\otimes s): a\in A; s\in S\}\subseteq A\otimes S$. Moreover, $\epsilon$ is
said to be {\it non-degenerate} (as a coaction) if the linear span of $N$ is
dense in $A\otimes S$.

\noindent (c) $v\in M(A\otimes S)$ is said to be a {\it unitary
corepresentation of $S$} if $({\rm id}\otimes \delta_S)(v) = v_{12}v_{13}$.

\medskip

Throughout this paper, $S$ is a saturated Hopf $C^*$-algebra with coaction $%
\delta_S$. Note that any Hopf $C^*$-algebra defined by regular or manageable
multiplicative unitary is saturated. We first observe the following
consequence of this extra condition. This fact will be used several times in
this paper.

\medskip

\noindent {\bf Remark 1.2:} Let $N$ be a $S$-invariant subspace of $S^*$.
Fix $f\in N$ and $t\in S$ such that $f(t)=1$. Then for any given $s\in S$,
there exist $r_i$, $s_i\in S$ such that $\sum_i\delta_S(s_i)(1\otimes r_i)$
converges to $s\otimes t$. Hence $\sum_i({\rm id}\otimes r_i\cdot
f)\delta_S(s_i)$ converges to $s$ and the subspace $N\cdot S$ generated by
elements of the form $({\rm id}\otimes g)\delta_S(r)$ ($g\in N$; $r\in S$)
is dense in $S$. Similarly, $S\cdot N$ is dense in $S$.

\medskip

\noindent {\bf Lemma 1.3:} Any left (or right) coidentity of $S$ is a
two-sided coidentity. Moreover, any coidentity of $S$ is positive.

\noindent {\bf Proof:} The first statement follows from the fact that any $%
h\in S^*$ such that $h S^*=(0)$ is zero (by Remark 1.2). For the second
statement, let $u\in S^*$ be a coidentity of $S$. Take any $f\in S^*_+$ and
any $a\in S_+$, let $b\in S$ be such that $a=b^*b$. By a similar argument as
in Remark 1.2, $a$ can be approximated by elements of the form $\sum_{ij}(%
{\rm id}\otimes r_j\cdot f\cdot r_i^*)\delta_S(c_i^*c_j)$. Therefore, $u(a)$
can be approximated by $\sum_{ij}(u\otimes r_j\cdot f\cdot r_i^*)\delta_S
(c_i^*c_j) = f((\sum_i c_ir_i)^*(\sum_j c_jr_j))\geq 0$.

\medskip

Note that $S^*$ is a Banach algebra with multiplication defined by $%
fg=(f\otimes g)\delta_S$ ($f,g\in S^*$). The argument of the above lemma
shows that if $N$ is any $S$-invariant subspace of $S^*$ and $u\in S^*$ such
that $uf = f$ for all $f\in N$, then $u\geq 0$ (note that $N$ is the 
vector subspace generated by $SN_+$).

\medskip

\noindent {\bf Lemma 1.4:} Let $A$ be a $S$-invariant subalgebra of $S^*$.
Any bounded left (or right) $\sigma(S^*,S)$-approximate identity of $A$ is a
bounded 2-sided $\sigma(S^*,S)$-approximate identity of $A$.

\noindent {\bf Proof:} Suppose that $\{f_i\}$ is a bounded left $%
\sigma(S^*,S)$-approximate identity of $A$. Since any $s\in S$ can be
approximated by sums of elements of the form $({\rm id}\otimes g)\delta_S(r)$
($g\in A$ and $r\in S$) (by Remark 1.2) and $\{h f_i\}$ is bounded for a
fixed $h\in A$, it is easy to see that $\{h f_i\}$ $\sigma(S^*,S)$-converges
to $h$.

\medskip

Now we will recall from [12] some basic materials about Fourier dualities of
Hopf $C^*$-algebras. We denote by $(\hat S_p, V_S)$ the {\it strong dual
object} of $S$ (i.e. the universal object corresponds to the unitary
corepresentations of $S$; see [12, 2.2] for a concrete definition). In this
case, $(\hat S_p, V_S, S)$ is called an {\it intrinsic dual pair}. We also
recall from [12, 2.7] the notion of {\it dualizable} Hopf $C^*$-algebra
which basically means that $\hat S_p$ is a Hopf $C^*$-algebra (see [12,
2.6]). Furthermore, $S$ is said to be {\it symmetrically dualizable} if the
intrinsic dual pair is a Fourier duality in the following sense. (Note that
if $S$ comes from a regular or manageable multiplicative unitary or their
tensor product, then $S$ is symmetrically dualizable.)

\medskip

\noindent {\bf Definition 1.5:} Let $S$ and $T$ be two Hopf $C^*$-algebras
and $X\in M(T\otimes S)$.

\noindent (a) $(T,X,S)$ is said to be a {\it Fourier duality} \vspace{-3mm}if

\begin{enumerate}
\item[i.]  $X$ and $X^\sigma\in M(S\otimes T)$ (where $\sigma$ is the flip
of the two variables) are unitary corepresentations of $S$ and $T$
respectively\vspace{-3mm};

\item[ii.]  $\{ ({\rm id}\otimes f)(X): f\in S^*\}\cap T$ and $\{(g\otimes 
{\rm id})(X): g\in T^*\}\cap S$ are dense in $T$ and $S$ respectively\vspace{%
-3mm}.
\end{enumerate}

Moreover, any $X\in M(T\otimes S)$ satisfying the first condition is called
a {\it unitary $T$-$S$-birepresentation}.

\noindent (b) Let $\mu$ and $\nu$ be representations of $S$ and $T$
respectively on $H$. Then $(\mu,\nu)$ is said to be a {\it $X$-covariant
representation} if $({\rm id}\otimes\mu)(X)_{12}X_{13}(\nu\otimes {\rm id}%
)(X)_{23} = (\nu\otimes {\rm id})(X)_{23}({\rm id}\otimes\mu)(X)_{12}$.

\noindent (c) A Fourier duality $(T,X,S)$ is said to be a Kac-Fourier
duality if there exists a $X$-covariant representation $(\mu,\nu)$ as well
as a $X^\sigma$-covariant representation $(\rho,\lambda)$ such that both $\mu
$ and $\nu$ are injective.

\medskip

In fact, the $X^\sigma$-covariant representation $(\rho,\lambda)$ in part
(c) is also injective (by Theorem 1.8 below).

\medskip

Note that we can define $V_S$-covariant representation in a similar way as
above even in the case when $S$ is not symmetrically dualizable.

\medskip

Now, if $(\mu,\nu)$ is a $X$-covariant representation, then by [12, 4.3], $%
V=(\nu\otimes\mu)(X)$ is a multiplicative unitary, i.e. $%
V_{12}V_{13}V_{23}=V_{23}V_{12}\in {\cal L}(H\otimes H\otimes H)$ (which
need not be semi-regular nor manageable) and $\mu(S)$ and $\nu(T)$ are Hopf $%
C^*$-algebras (which need not be defined by $V$ in the usual sense as in
[2]). Conversely, any good (in particular, regular or manageable)
multiplicative unitary defines a Fourier duality (which has a canonical
covariant representation) in the obvious way.

\medskip

Moreover, a Kac-Fourier duality is exactly a duality $(T,X,S)$ defined (in
the usual sense) by two multiplicative unitaries $V$ and $W$ such that $%
S_V\cong S\cong \hat S_W$. Hence Kac-Fourier duality is much more general
than Kac system. For the other properties of Fourier duality, we refer the
readers to [12].

\medskip

From now on, $(T,X,S)$ is a Fourier duality. Let $S^\#$ be the set $\{f\in
S^*: ({\rm id}\otimes f)(X)\in T\}$ and $T^\#=\{g\in T^*:(g\otimes{\rm id}%
)(X)\in S\}$. Moreover, we denote by $j_S$ and $j_T$ the canonical maps from 
$T^\#$ and $S^\#$ to $S$ and $T$ respectively.

\medskip

\noindent {\bf Lemma 1.6:} For a Fourier duality $(T,X,S)$, any coidentity $u
$ of $S$ is a Hopf $*$-homomorphism and $({\rm id}\otimes u)(X) = 1$.

\noindent {\bf Proof:} Since $u$ is a coidentity, $(u\otimes {\rm id}%
)\delta_S = {\rm id}$. Hence $X = ({\rm id}\otimes {\rm id}\otimes u)({\rm id%
} \otimes \delta_S)(X) = X(({\rm id}\otimes u)(X)\otimes 1)$ and so $({\rm id%
}\otimes u)(X) = 1$. Moreover, for any $f,g\in T^*$, $u(j_S(f)j_S(g)) = (f
g\otimes u)(X) = f g(1) = f(1)g(1) = u(j_S(f))u(j_S(g))$.

\medskip

Hence if $S$ is symmetrically dualizable, then any coidentity of $S$ is a
Hopf $*$-homomorphism. Next, we would like to recall the notion of Fourier
algebras. A non-zero $S$-invariant 2-sided (respectively, left or right)
ideal of $S^*$ is called a {\it 2-sided} (respectively, {\it left} or {\it %
right}) {\it Hopf ideal} of $S^*$.

\medskip

\noindent {\bf Definition 1.7:} The intersection, $A_S$ (respectively, $A_S^l
$ and $A_S^r$), of all 2-sided (respectively, left and right) Hopf ideals of 
$S^*$ is called the {\it Fourier algebra} (respectively, {\it left} and {\it %
right Fourier algebra}) of $S$.

\medskip

We recall from [12, 5.5 \& 1.9] that if $(\mu,\nu)$ is a $X$-covariant
representation, then $\nu^*({\cal L}(H)_*)$ and $\mu^*({\cal L}(H)_*)$ are
left Hopf ideal of $T^*$ and right Hopf ideal of $S^*$ respectively.
Moreover, we recall the following result from [12, 5.9].

\medskip

\noindent {\bf Theorem 1.8:} Let $(T,X,S)$ be a Kac-Fourier duality and let $%
(\mu,\nu)$ and $(\rho,\lambda)$ be a $X$-covariant and a $X^\sigma$%
-covariant representations respectively. Then \vspace{-2.5mm} 
\[
e_\mu =e_{\lambda} \qquad {\rm and} \qquad A^l_S=A^r_S=A_S=e_\mu S^*=\mu^*(%
{\cal L}(H)_*)\vspace{-2.5mm}
\]
(where $e_\mu$ and $e_\lambda$ are the support projections of $\mu$ and $%
\lambda$ respectively). The same is true for $T$.

\medskip

Note that we can deduce the unique of Kac-Fourier duality as well as many
other interesting properties from this theorem (see [12, \S 5]).

\medskip

We end this section with some notations and results about Banach algebras
which will be used in the following sections. From now on, until the end of this
section, $A$ is a Banach algebra.

\medskip

\noindent {\bf Definition 1.9:} Let $M_l(A) = \{l\in {\cal L}(A): l(ab)=
l(a)b; a,b\in A\}$ (i.e. the set of all {\it left multipliers} of $A$). Let $%
{\cal S}_l$ be the topology on $M_l(A)$ corresponding to the family of
semi-norms $\{p_a:a\in A\}$ (where $p_a(l) = \|l(a)\|$ for any $l\in M_l(A)$%
). Then ${\cal S}_l$ is called the {\it left strict topology} on $M_l(A)$.

\medskip

Suppose that $A$ is a left ideal of a Banach algebra $B$. Then there exists
a canonical map from $B$ to $M_l(A)$. Moreover, we have the following well
know result.

\medskip

\noindent {\bf Lemma 1.10:} If $A$ has a left approximate identity, then the
canonical image of $A$ is dense in $(M_l(A), {\cal S}_l)$.

\medskip

We can also define $M_r(A)$ as well as $M(A)$ (i.e. $M(A;A)$ as defined in
Definition 1.11 below) and have similar density results for them. More
generally, for any Banach $A$-bimodules $X$ and $Y$, we denote by $M^r_A(X,Y)
$ (respectively, $M^l_A(X,Y)$) the set of all elements in ${\cal L}(X,Y)$
(i.e. bounded linear operators from $X$ to $Y$) that respect the left
(respectively, right) $A$-multiplication.

\medskip

\noindent {\bf Definition 1.11:} Let $X$ be a Banach $A$-bimodule. Let $l\in
M^l_A(A,X)$ and $r\in M^r_A(A,X)$. Then $(l,r)$ is said to be a {\it %
multiplier} if $a\cdot l(b) = r(a)\cdot b$ for any $a,b\in A$. We denote by $%
M(A,X)$ the set of all multipliers of $X$.

\medskip

It is clear that all of the spaces $M^r_A$, $M^l_A$ and $M$ are Banach
spaces (under the canonical norms on them). Moreover, $M^r_A(X,Y^*)\cong
M^l_A(Y,X^*)$. We will now give some basic facts about approximate
identities and invariant means.

\medskip

\noindent {\bf Proposition 1.12:} (a) The following statements are
equivalent:

\noindent (i) $A$ has a bounded right approximate identity;

\noindent (ii) $A$ has a bounded right weak approximate identity;

\noindent (iii) $(A^{**}, \times_1)$ has a right identity (where $\times_1$
is the first Arens multiplication on $A^{**}$);

\noindent (iv) for any Banach $A$-bimodule $X$, $M^r_A(A, X^*)\cong X^*$
canonically;

\noindent (v) $A^{**} \cong M^r_A(A,A^{**})$ canonically;

\noindent (b) If $A$ has both a bounded right approximate identity and a
bounded left approximate identity, then $A$ has a bounded two-sided
approximate identity.

\noindent {\bf Proof:} Since part (b) is well known, we need only to show
part (a). In fact, the equivalences of (i)-(iii) in part (a) are also well
known (see e.g. [13, 5.1.8]). To prove that (i) implies (iv), we let $\{a_i\}
$ be a bounded right approximate identity of $A$ and $r\in M^r_A(A, X^*)$.
Then $\{ r(a_i)\}$ is a bounded net in $X^*$ and hence has a subnet
converging weakly to an element $\varphi_0$ of $X^*$. It is easy to check
that $r(a) = a\cdot \varphi_0$ and (iv) holds. Condition (iv) obviously
implies condition (v). Finally, we show that (v) implies (iii). If $u$ is an
element in $A^{**}$ corresponding to ${\rm id}_A\in M^r_A(A,A^{**})$, then $u
$ is clearly a right identity for the first Arens product.

\medskip

A cone $E_+$ of a normed space $E$ is said to be a $L$-cone if $\|x+y\| =
\|x\| + \|y\|$ for any $x,y \in E_+$.

\medskip

\noindent {\bf Definition 1.13:} (a) A Banach algebra $A$ with a closed cone 
$A_+$ is said to be a {\it $LB$-algebra} if $A_+$ is a generating $L$-cone
and there exists $e\in A^*_+$ with $e(v)=\|v\|$ for all $v\in A_+$.

\noindent (b) Let $A$ be a $LB$-algebra and $S=\{x\in A_+: \|x\|=1\}$.

\noindent (i) A net $\{a_i\}$ in $S$ is said to be a {\it left}
(respectively, {\it left weak}) {\it approximate invariant mean} if for 
any $a\in S$, $aa_i - a_i$ converges to 
$0$ in norm (respectively, in the weak topology). We can
also define a {\it right {\rm ( respectively,} right weak{\rm )} approximate
invariant mean} in a similar way.

\noindent (ii) $m\in A^{**}_+$ is said to be a {\it left} (respectively, 
{\it right}) {\it invariant mean} if $\|m\|=1$ and $a\cdot m = e(a)m$
(respectively, $m\cdot a = e(a)m$) for any $a\in A$.

\medskip

Note that the predual of a Hopf von Neumann algebra is a $LB$-algebra. If $G$
is a locally compact group and $A=L^1(G)$, then $m\in A^{**}=L^\infty(G)^*$
is a left invariant mean in the above sense if and only if it is a left
invariant mean in the usual sense.

\medskip

\noindent {\bf Proposition 1.14:} Let $A$ be a $LB$-algebra.

\noindent (a) The following statements are equivalent.

\hspace{-1em} (i) $A^{**}$ has a left (respectively, right) invariant mean.

\hspace{-1em} (ii) $A$ has a left (respectively, right) weak approximate
invariant mean.

\hspace{-1em} (iii) $A$ has a left (respectively, right) approximate
invariant mean.

\noindent (b) Suppose that $A$ satisfies the extra assumption that $%
\|ab\|=\|a\|\,\|b\|$ for any $a,b\in A_+$. Then $A$ has a 2-sided
(respectively, weak) approximate invariant mean if and only if $A$ has both
a left and a right (respectively, weak) approximate invariant means.

\noindent {\bf Proof:} (a) If $m$ is a left invariant mean of $A^{**}$, then
there exist $b_i\in A_+$ such that $\|b_i\|\leq 1$ and $\{b_i\}$ $%
\sigma(A^{**},A^*)$-converges to $m$. Hence $\|b_i\| = e(b_i)$ converges to $%
m(e) =1$ and we can assume all $b_i$'s to be non-zero. It is easily seen
that $a_i=b_i/\|b_i\|$ is a left weak approximate invariant mean. The proof
for (ii) imply (iii) follows from a similar argument as in [5, 2.8.4].
Finally, (iii) implies (i) can be shown easily using the weak compactness of
the unit ball of $A^{**}$.

\noindent (b) Let $\{a_i\}$ and $\{b_j\}$ be respectively a left and a right
approximate invariant means. Then clearly $\{a_ib_j\}$ is a 2-sided
approximate invariant mean.

\medskip

\medskip

\medskip

\noindent {\large 2. Coamenability of Hopf $C^*$-algebras}

\medskip

\medskip

In this section, we will study coamenability of Hopf $C^*$-algebras. We will
give three definitions for it, namely, $H_1$-coamenable, $H_2$-coamenable and $%
H_3$-amenable. We will then study the relations between them. In particular,
we will show that $H_1$-coamenability is stronger than both $H_2$-coamenability
and $H_3$-coamenability. Moreover, 2-sided $H_2$-coamenability will imply
2-sided $H_3$-coamenability. We will also study some properties and
implications of these coamenabilities.

\medskip

\noindent {\bf Definition 2.1:} A (saturated) Hopf $C^*$-algebra $S$ is said
to be:

\noindent (a) {\it 2-sided coamenable as Hopf $C^*$-algebra of the first type}%
, or simply {\it 2-sided $H_1$-coamenable}, if $A_S$ is non-zero and has a
bounded 2-sided approximate identity.

\noindent (b) {\it 2-sided coamenable as Hopf $C^*$-algebra of the second type}%
, or simply {\it 2-sided $H_2$-coamenable}, if $S^*$ is unital and any 2-sided
Hopf ideal separates points of $S$.

\noindent (c) {\it 2-sided coamenable as Hopf $C^*$-algebra of the third type}%
, or simply {\it 2-sided $H_3$-coamenable}, if the canonical map from $S^*$ to 
$M(A_S)$ is surjective.

\medskip

\noindent {\bf Remark 2.2} We can also define 1-sided version of Definition
2.1: $S$ is said to be:

\noindent (i) left $H_1$-coamenable if $A_S^l$ has a bounded left approximate
identity.

\noindent (ii) left $H_2$-coamenable if $S^*$ is unital and any left Hopf
ideal separates points of $S$.

\noindent (iii) left $H_3$-coamenable if the canonical map from $S^*$ to $%
M_l(A_S^l)$ is injective.

\noindent Similarly, we can define right $H_i$-coamenability ($i=1,2,3$).

\medskip

Note that since $A_S\cdot S$, $A_S^l\cdot S$ and $S\cdot A_S^r$ are dense in 
$S$ (by Remark 1.2), any $f\in S^*$ such that $fA_S=0$ or $fA_S^l=0$ or $%
A_S^rf=0$ is zero. Hence the canonical map in Definition 2.1(c) as well as
its one-sided companions are automatically injective.

\medskip

\noindent {\bf Example 2.3:} (a) Let $G$ be a locally compact group. If $%
S=C_0(G)$, then $S^*=M(G)$ and $A_S = A_S^l = A_S^r = L^1(G)$ (by Theorem
1.8). It is clear that $S$ is automatically 2-sided (and also left and
right) $H_i$-coamenable for $i=1,2,3$. On the other hand, let $S=C^*(G)$ (or $%
S=C^*_r(G)$). Then $S^*=B(G)$ by definition (respectively, $S^*=B_r(G)$) and 
$A_S = A_S^l = A_S^r = A(G)$ (by Theorem 1.8 or [4]). It is well known that
for any fixed $i=1,2,3$, $S$ is 2-sided $H_i$-coamenable if and only if $G$ is
amenable (see e.g. [9] and [8, 4.1]).

\noindent (b) Let $S_V$ and $\hat S_V$ be the Hopf $C^*$-algebras defined by
a biregular irreducible multiplicative unitary $V$. Then $\hat S_V$ is 
$H_i$-coamenable ($i=1,2,3$) if and only if $V$ is amenable in the sense of [2,
A13(c)] (by Theorem 1.8, Proposition 3.5 and Corollary 3.13 below; 
see also Examples 3.2(c) as well as the discussion after Definition 1.5). In
particular, any discrete quantum group (i.e. the dual Hopf $C^*$-algebra 
of a unital Hopf $C^*$-algebra; see [19]) is $H_i$-coamenable ($i=1,2,3$).

\medskip

These three coamenabilities look very distinct. However, we will see in the
following (Theorem 2.8) some relations between them. We will also show later
that in the case of Kac-Fourier duality, they are essentially the same (see
Corollary 3.13 and Remark 3.15 below).

\medskip

\noindent {\bf Proposition 2.4:} Let $A$ be a left Hopf ideal of $S^*$. Then 
$S^*$ is unital and $\overline{A}^{\sigma(S^*,S)} =S^*$ if and only if $A$
has a bounded left $\sigma(S^*,S)\!\mid_A$-approximate identity.

\noindent {\bf Proof:} Suppose that $\{f_i\}\subseteq A$ is a bounded left $%
\sigma(S^*,S)$-approximate identity. Then there exists a subnet $\{g_j\}$ of 
$\{f_i\}$ such that $\{g_j\}$ is $\sigma(S^*,S)$-convergent to an element $%
g_0\in S^*$. For any $f\in A$, $s\in S$ and $h\in S^*$, \vspace{-2mm}
\begin{eqnarray*}
(g_0 h)(f\cdot s) & = & \lim_j g_j(h\cdot(f\cdot s)) \; = \; \lim_j (g_j(h
f))(s)\vspace{-1mm} \\
& = & (h f)(s) \; = \; h(f\cdot s). \vspace{-3mm}
\end{eqnarray*}
Hence, by Remark 1.2, $g_0$ is a left identity of $S^*$ and so is a 2-sided
identity by Lemma 1.3. Moreover, by a similar argument as above, it is not
hard to see that for any $h\in S^*$, $hf_i\in A$ will $\sigma(S^*,S)$%
-converge to $h$. Conversely, suppose that $S^*$ is unital and $A$ is $%
\sigma(S^*,S)$-dense in $S^*$. As $A$ separates points of $S$, the canonical 
$*$-homomorphism $i_S$ from $S$ to $A^*$ is injective (and so a metric
injection). Hence the map $i_S^*$ from $A^{**}$ to $S^*$ is a metric
surjection and is $\sigma(A^{**},A^*)$-$\sigma(S^*,S)$ continuous. Let $U$
and $U_A$ be respectively the closed unit balls of $S^*$ and $A$. Then $%
U_A=i_S^*(\kappa_A(U_A))$ is $\sigma(S^*,S)$-dense in $U$ (where $\kappa_A$
is the canonical map from $A$ to $A^{**}$). Therefore, if $f_0$ is the
identity of $S^*$, there exists a net $\{f_i\}$ in $A$ such that $%
\|f_i\|\leq\|f_0\|$ and $\{f_i\}$ $\sigma(S^*,S)$-converges to $f_0$. It is
clear that $f_i$ is a bounded left $\sigma(S^*,S)$-approximate identity for $%
A$ (in fact, for $S^*$).

\medskip

We have a similar result for right (as well as 2-sided) Hopf ideals of $S^*$%
. Note that the topology considered in Proposition 2.4 (i.e. $%
\sigma(S^*,S)\!\mid_A$) is strictly weaker than $\sigma(A,A^*)$ (consider $S$
as a subset of $A^*$). Therefore, (1-sided and 2-sided) $H_1$-coamenability is
stronger than $H_2$-coamenability. However, there is no reason to believe that
the existence a bounded $\sigma(S^*,S)\!\mid_A$-approximate identity will
imply the existence of a bounded approximate identity. Therefore, 
$H_1$-coamenability and $H_2$-coamenability are probably distinct in general.

\medskip

On the other hand, Proposition 2.4 and Lemma 1.4 also imply the following
unexpected result (note that a right approximate identity is considered
instead of a left one).

\medskip

\noindent {\bf Corollary 2.5:} Let $A$ be a left Hopf ideal of $S^*$ with a
bounded right $\sigma(S^*,S)$-approximate identity (in particular, if $A$
has a bounded right approximate identity). Then $S^*$ equals $\overline{A}%
^{\sigma(S^*,S)}$ and is unital.

\medskip

Thus, in the case when $A=S^*$, we obtain the following interesting
proposition.

\medskip

\noindent {\bf Proposition 2.6:} Let $S$ be a (saturated) Hopf $C^*$%
-algebra. The following statements are equivalent.

\noindent (i) $S^*$ is unital.

\noindent (ii) $S^*$ has a bounded left (or right) $\sigma(S^*,S^{**})$%
-approximate identity.

\noindent (iii) $S^*$ has a bounded left (or right) approximate identity.

\noindent (iv) $S^*$ has a bounded left (or right) $\sigma(S^*,S)$%
-approximate identity.

\medskip

In fact, it is clear that (i) $\Rightarrow$ (iii) $\Rightarrow$ (ii) $%
\Rightarrow$ (iv) and Proposition 2.4 shows that (iv) implies (i).

\medskip

Next, we will study the relation between $H_2$-coamenability and 
$H_3$-coamenability.

\medskip

\noindent {\bf Proposition 2.7:} Let $A$ be a 2-sided Hopf ideal of $S^*$.
Suppose that $S^*$ is unital and $A$ is $\sigma(S^*,S)$-dense in $S^*$. Then 
$M(A)\cong S^*$ canonically.

\noindent {\bf Proof:} It is required to show that the canonical map from $%
S^*$ to $M(A)$ is surjective. Let $(l,r)$ be an element in $M(A)$. We recall
from Remark 1.2 that any $s\in S$ can be approximated by sums of elements of
the form $(f\otimes {\rm id})\delta_S(t)$ ($f\in A$, $t\in S$). As the
canonical map $j$ from $S$ to $A^*$ (which is the composition of $\kappa_S$
with the canonical map from $S^{**}$ to $A^*$) respects the
comultiplications (i.e. $\delta_{A^*}\circ j = (j\otimes j)\circ\delta_S$),
we see that $j(s)$ can be approximated (in norm) by sums of elements of the
form $(f\otimes {\rm id}) \delta_{A^*}(j(t))$. Moreover since for any $%
f,g\in A$ and $t\in S$, $l^*(j(t)\cdot f)(g) = (fl(g))(j(t)) = (r(f)g)(j(t))
= (g\cdot j(t))(r(f))$, \vspace{-1.5mm} 
\[
l^*((f\otimes {\rm id})\delta_{A^*}(j(t))) = (r(f)\otimes {\rm id}%
)\delta_{A^*}(j(t))\vspace {-1.5mm}
\]
which is in $j(S)$. Thus, $l^*(j(S))\subseteq j(S)$. Similarly, we have $%
r^*(j(S))\subseteq j(S)$. Therefore, as $j$ is injective (for $A$ separates
points of $S$), we can regard $l^*, r^*\in {\cal L}(S)$. Note that under
this identification, we still have $g\cdot l^*(s) = l^*(g\cdot s)$ (as $j$
respects the comultiplication and $A$ separates points of $S$). Now, if $u$
is the identity in $S^*$, then $h=u\circ l^*\in S^*$. For any $g\in A$ and $%
s\in S$\vspace{-1mm}, 
\[
g(l^*(s)) = (u\otimes g)\delta_S(l^*(s)) = u(g\cdot l^*(s)) = h(g\cdot s).%
\vspace{-0.5mm}
\]
Therefore, $l(g) = h g$. On the other hand, let $\{ g_i\}$ be a net in $A$
that $\sigma(S^*,S)$-converges to $u$. For any $f\in A$, we have $l^*(s\cdot
f)(g_i) = r^*(g_i\cdot s)(f)$. The left hand side of the equality converges
to $h(s\cdot f)$ while the right hand side equals $(r(f)\otimes
g_i)\delta_S(s)$ which converges to $r(f)(s)$. Hence $r(f) = f h$. Now the
image of $h\in S^*$ is $(l,r)\in M(A)$ and the result is proved.

\medskip

It is not known whether Proposition 2.7 holds for left (or right) Hopf
ideals in general. Note that it is not clear whether the image of $\overline{%
A}^{\sigma(S^*,S)}$ is a left ideal of $M_l(A)$. Nevertheless, we can still
prove that one-sided $H_1$-coamenability is stronger that $H_3$-coamenability.

\medskip

\noindent {\bf Theorem 2.8} (a) Left (respectively, right or 2-sided) 
$H_1$-coamenability will imply left (respectively, right or 2-sided) 
$H_2$-coamenability.

\noindent (b) 2-sided $H_2$-coamenability implies 2-sided $H_3$-coamenability.

\noindent (c) Left (or right) $H_1$-coamenability will imply left
(respectively, right) $H_3$-coamenability.

\noindent {\bf Proof:} Because of Proposition 2.7 and the paragraph after
Proposition 2.4, we need only to prove part (c). Let $A$ be a left Hopf
ideal of $S^*$ with a bounded left approximate identity $\{f_i\}$. Note that 
$A$ can be regarded as a subspace of $M_l(A)$ (as the canonical map is
injective). By Remark 1.2, $A\cdot S$ is dense in $S$. Hence $\{f_i\}$ is an
approximate identity for the Banach left $A$-module $S$. Suppose that $%
\{g_j\}$ is a net in $A$ that ${\cal S}_l$-converges to an element $g\in A$
(see Definition 1.9). For any $s\in S$, since $s=h\cdot t$ for some $h\in A$
and $t\in S$, $g_j(s) = g_jh(t)$ converges to $gh(t) = g(s)$. This shows
that $\sigma(S^*,S)\!\mid_{A}$ is weaker than ${\cal S}_l\!\mid_{A}$. Hence
Lemma 1.10 ensures that there is a continuous map $\Psi$ from $M_l(A)$ to
the completion of $(A,\sigma(S^*,S)\!\mid_A)$ which can be regarded as a
subspace of the algebraic dual $S^d$ of $S$. Now for any $l\in M_l(A)$, $%
g_i=l(f_i)\in A$ will ${\cal S}_l$-converge to $l$ and so $\Psi(l)(s) =
\lim_i g_i(s)$ for any $s\in S$. Therefore, $\|\Psi(l)(s)\|\leq
\|l\|\sup_i\|f_i\|\|s\|$ which implies that $\Psi(l)\in S^*$. Thus, $\Psi$
is a continuous map from $(M_l(A), {\cal S}_l)$ to $(S^*, \sigma(S^*,S))$
which extends the identity map on $A$. Now let $j$ be the canonical
injection from $S^*$ to $M_l(A)$. For any $l\in M_l(A)$ and $g\in A$\vspace{%
-1.5mm}, 
\[
j(\Psi(l))(g) = \Psi(l)g = w-\lim_i l(f_i)g = w-\lim_i l(f_ig).\vspace
{-2.5mm}
\]
However, since $f_ig$ converges to $g$, $l(f_ig)$ will converge to $l(g)$ in
norm. Therefore, $j\circ \Psi= {\rm id}$ and $j$ is surjective and hence
bijective. Now, part (c) is proved by putting $A=A_S^l$.

\medskip

We now consider some implications of these coamenabilities. The following
lemma is clear.

\medskip

\noindent {\bf Lemma 2.9:} Let $A$ be a left Hopf ideal of $S^*$. Then each
of the following conditions is weaker than the previous ones.

\noindent (i) $S^*= M_l(A)$.

\noindent (ii) the norm given by $\|a\|_l = \sup\{\|ab\|:b\in A; \|b\|\leq
1\}$ is equivalent to the original norm on $A$.

\noindent (iii) $A$ is norm closed in $(M_l(A), \|\cdot\|_l)$.

\medskip

Note that in the case when $S=C^*(G)$ and $A=A(G)$, the above conditions are
all equivalent to the amenability of $G$ (see [9]). However, we do not know
whether it is the case in this general situation.

\medskip

\noindent {\bf Lemma 2.10:} Let $(T,X,S)$ be a Fourier duality such that $%
j_S(A_T^r)$ (see Section 1) is dense in $S$. If $T$ is right $H_1$-coamenable,
then the following conditions hold.

\noindent (i) $A_T^r$ has a bounded approximate identity for the right $A_T^r
$-module $S$ (the module structure is given by $j_S(g)\bullet h = j_S(g h)$
for any $g\in T^\#$ and $h\in A_T^r$).

\noindent (ii) If $N$ is any $S$-invariant subspace of $S^*$, then $%
M^r_{A_T^r}(A_T^r, N) \cong N$.

\noindent {\bf Proof:} (i) We first note that by the definition of Fourier
duality (see Definition 1.5(a)), $j_S(T^\#)$ is dense in $S$. It is clear
that $\|j_S(g)\bullet h\| \leq \|j_S(g)\|\|j_S(h)\| \leq \|j_S(g)\|\|h\|$.
Since $A_T^r$ has a bounded right approximate identity and $j_S(T^\#)\bullet
A_T^r = j_S(T^\# A_T^r)$, $A_T^r A_T^r = A_T^r$ and $S\bullet A_T^r$ is
dense in $S$ (as $j_S(A^r_T)$ is dense in $S$ by assumption). Therefore, the
bounded right approximate identity of $A_T^r$ is a bounded approximate
identity for $S$.

\noindent (ii) Note that as $j_S(A^r_T)$ is dense in $S$, the corresponding
dual left $A_T^r$-module structure on $S^*$ is given by $h\bullet f =
j_S(h)\cdot f\in S^*$ ($h\in A_T^r$; $f\in S^*$). Thus, $N$ is a $A_T^r$%
-submodule of $S^*$. By Proposition 1.12(a)(iv), $M^r_{A_T^r}(A_T^r,
S^*)\cong S^*$. Hence for any $r\in M^r_{A_T^r}(A_T^r, N)$, there exists $%
f_0\in S^*$ such that $r(h) = h\bullet f_0 \in N$ for any $h\in A_T^r$.
Since $j_S(A_T^r)$ is dense in $S$, $S\cdot f_0 \subseteq N$ which implies
that $f_0\in N$.

\medskip

In the case when $S=C_0(G)$, $A_T^r=A(G)$ and $N=L^1(G)$, the above two
conditions are also equivalent to the amenability of $G$ (see [10, Theorem
1]). Again, we do not know whether it is true in general.

\medskip

\medskip

\medskip

\noindent {\large 3. Amenable covariant representations}

\medskip

\medskip

In this section, we will define and study amenable covariant
representations. In particular, we will give a partial analogue of [5, 2.4].
Moreover, we will study the crossed products as in [11] and show that
nuclearity is preserved under crossed products of coactions of amenable
Hopf $C^*$-algebras. We will then apply these results to Kac-Fourier
dualities and show that the three coamenabilities in Section 2 are 
essentially the same in this case (see Corollary 3.13 and Remark 3.15).

\medskip

Throughout this section, unless specified, $S$ is a (saturated) dualizable
Hopf $C^*$-algebra (which is the case if $S$ comes from a regular or
manageable multiplicative unitary) and $(\hat S_p, V_S, S)$ is the intrinsic
duality of $S$ (see Section 1).

\medskip

\noindent {\bf Definition 3.1:} (a) A $V_S$-covariant representation $(\mu,
\nu)$ (see the paragraph after Definition 1.5) is said to be {\it amenable}
if $\nu$ is injective.

\noindent (b) $S$ is said to be {\it amenable} if there exists an amenable 
$V_S$-covariant representation.

\medskip

\noindent {\bf Example 3.2} (a) $C_r^*(G)$ and $C^*(G)$ are automatically
amenable. In fact, the unitary corepresentations of $C_r^*(G)$ and $C^*(G)$
are the same. Let $(C_0(G), V_G)$ and $(C_0(G), V_G^r)$ be the strong dual
object of $C^*(G)$ and $C^*_r(G)$ respectively. Since $V_G^r$ is a quotient
of $V_G$ and there exists an amenable $V_G^r$-covariant representation, $%
C^*(G)$ and $C^*_r(G)$ are amenable.

\noindent (b) Let $(C^*(G), W_G)$ be the strong dual of $C_0(G)$. Theorem
1.8 implies that for any $W_G$-covariant representation $(\mu,\nu)$, $%
\nu(C^*(G)) = C_r^*(G)$. Therefore, $C_0(G)$ is amenable if and only if $G$
is amenable (using Theorem 1.8 again).

\noindent (c) Let $V$ be a regular or manageable multiplicative unitary and $%
S = S_V$. If $V$ is amenable in the sense of [2, A13(c)], then $S_V$ is
amenable. If, in addition, $V$ is irreducible, then the amenability of $%
S_V$ will imply the amenability of $V$ (by Theorem 1.8).

\medskip

\noindent {\bf Lemma 3.3:} Suppose that $(\mu,\nu)$ is a $V_S$-covariant
representation of $S$ on a Hilbert space $H$. Let $T_0=\nu(\hat S_p)$ and $B
= \nu^*({\cal L}(H)_*)$. Then $\nu^*(T_0^*)$ and $B$ are both left Hopf
ideals of $\hat S_p^*$ such that $\nu^*(T_0^*)= \overline{B}^{\sigma(\hat
S_p^*,\hat S_p)}$.

\noindent {\bf Proof:} By [12, 5.5 \& 1.9], $B$ is a left Hopf ideal of $%
\hat S_p^*$. It is not hard to see that $\nu^*(T_0^*)= \overline{B}%
^{\sigma(\hat S_p^*,\hat S_p)}$ and hence $\nu^*(T_0^*)$ is also a left Hopf
ideal of $\hat S_p^*$.

\medskip

\noindent {\bf Corollary 3.4:} If $\hat S_p$ is left $H_2$-coamenable, then
any $V_S$-covariant representations is amenable.

\medskip

Similar to the case of Kac algebras (see [5, \S 2]), we have the following
list of equivalent conditions for amenable $V_S$-covariant representations.
Note that these equivalences are known in the case of regular multiplicative
unitary (see [3, 5.5]).

\noindent

\medskip

\noindent {\bf Proposition 3.5:} Let $(\mu, \nu)$ be a $V_S$-covariant
representation on $H$ and $V=(\nu\otimes\mu)(V_S)$. Let $T_0=\nu(\hat S_p)$
and $B=\nu^*({\cal L}(H)_*)$. Then the following statements are equivalent. 
\vspace{-3mm}

\begin{enumerate}
\item[i.]  $(\mu,\nu)$ is amenable\vspace{-3mm}.

\item[ii.]  $T_0^*$ is unital\vspace{-3mm}.

\item[iii.]  The trivial representation $\pi_e$ of $\hat S_p$ is weakly
contained in $\nu$\vspace{-3mm}.

\item[iv.]  There exists a net of unit vectors $\{\eta_i\}$ in $H$ such that 
$\|V(\eta_i\otimes\xi) -\eta_i\otimes\xi\|$ converges to zero for any $%
\xi\in H$\vspace{-3mm}.

\item[v.]  $B$ has a bounded left approximate identity\vspace{-3mm}.
\end{enumerate}

\noindent {\bf Proof:} It is clear that condition (i) implies both (ii) and
(iii). To show that (ii) implies (i), we first note that, by the proof of
the equality in Lemma 1.6, the image of the identity of $T_0^*$ in $\hat
S_p^*$ is an identity. Now Lemma 3.3 shows that condition (i) holds. We can
use a similar argument as [5, 2.7.1] and [5, 2.7.2] to prove that condition
(iii) is stronger than (iv) and condition (iv) is stronger than (v).
Finally, it follows from Lemma 3.3 and Proposition 2.4 that condition (v)
implies condition (i).

\medskip

Note that the condition of $B$ having a bounded right approximate identity
is stronger than it having a bounded left approximate identity (see
Proposition 3.6(b) below).

\medskip

Suppose that $A=\mu^*({\cal L}(H)_*)$. Then $A^*\cong
\mu(S)^{\prime\prime}\subseteq {\cal L}(H)_*$ is a Hopf von Neumann algebra
and hence $A$ is a $LB$-algebra (see Definition 1.3(a)). We have a sort of
invariant mean characterisation for amenability.

\medskip

\noindent {\bf Proposition 3.6:} Let the situation be the same as
Proposition 3.5.

\noindent (a) Let $A=\mu^*({\cal L}(H)_*)$ and $j$ be the canonical
non-degenerate $*$-homomorphism from $S$ to $A^*$. Then conditions (i)-(v)
imply condition (vi) below. If, in addition, $\{({\rm id}\otimes\omega\circ%
\mu)(V_S):\omega\in {\cal L}(H)_*\}$ is dense in $\hat S_p$, then conditions
(i)-(vi) are equivalent to one another\vspace{-2mm}.

\begin{enumerate}
\item[vi.]  $A^{**}$ has a right invariant mean $m$ such that $%
m(\{(g\circ\nu\otimes j)(V_S): g\in T_0^*\}) \neq (0)$\vspace{-2mm}.
\end{enumerate}

\noindent (b) If $B=\nu^*({\cal L}(H)_*)$ has a bounded right 
approximate identity, then $(\mu,\nu)$ is amenable.

\noindent {\bf Proof:} (a) We will first show that (iv) implies (vi). For
any $\zeta\in H$, let $\mu_\zeta = \mu^*(\omega_\zeta)$. Since $%
\{\mu_{\eta_i}\}\subseteq A$ is a bounded net, it has a subnet $%
\{\mu_{\zeta_j}\}$ $\sigma(A^{**},A^*)$-converging to $m\in A^{**}$. For any 
$\xi\in H$ and $x\in A^*$, we have \vspace{-3mm}
\begin{eqnarray*}
(m\cdot\mu_\xi)(x) & = & \lim_j \omega_{\zeta_j\otimes\xi}
(V(\tilde\mu(x)\otimes 1)V^*) \; = \; \lim_j \langle (\tilde\mu(x)\otimes
1)V^*(\zeta_j\otimes \xi), V^*(\zeta_j\otimes\xi)\rangle \\
& = & \lim_j \langle (\tilde\mu(x)\otimes 1)(\zeta_j\otimes \xi),
\zeta_j\otimes \xi\rangle \; = \; \mu_\xi(1) m(x)\vspace {-2mm}
\end{eqnarray*}
(where $\tilde\mu$ is the induced map from $A^*$ to ${\cal L}(H)$ and 
$V=(\nu\otimes\mu)(V_S)$). Hence, $%
m\cdot a = a(1)m$ if $a$ is a linear combination of $\mu_\xi$'s. Therefore, $%
m$ is a right invariant mean for $A^{**}$ (as $m$ is continuous and $(m\cdot
a)(z) = m(a\cdot z)$). Moreover, as $m(1) = 1$ and $T_0^*$ is unital, we
have by the proof of the equality in Lemma 1.6, $m(\{(g\otimes j)(V): g\in
T_0^*\}) \neq (0)$. Conversely, suppose that condition (vi) is true. Let $%
\{a_i\}$ be a net in $A$ bounded by $K=\|m\|$ and $\sigma(A^{**},A^*)$%
-converging to $m$. Take any $f\in T_0^*$ such that $m((f\circ\nu\otimes
j)(V_S))=1$. Then for any $a\in A\subseteq S^*$, \vspace{-3mm}
\begin{eqnarray*}
\mid\!\pi_e(j_{\hat S_p}(a))\!\mid & = & \mid\! a(1)\!\mid\;\; =\;\; \mid\!
m\cdot a((f\circ\nu\otimes j)(V_S))\!\mid \\
& = & \mid\!\lim_i (f\circ\nu\otimes a_i\otimes a)
((V_S)_{12}(V_S)_{13})\!\mid \\
& = & \lim_i\mid\! f\circ\nu(j_{\hat S_p}(a_i)j_{\hat S_p}(a))\!\mid \quad
\leq \quad K\|f\|\,\|\nu(j_{\hat S_p}(a))\|
\end{eqnarray*}
(where $j_{\hat S_p}(a) = ({\rm id}\otimes a)(V_S)$). Since $j_{\hat S_p}(A)$
is dense in $\hat S_p$, this shows that condition (iii) holds.

\noindent (b) This follows directly from Corollary 2.5 and Lemma 3.3.

\medskip

\noindent {\bf Remark 3.7:} (a) Let $V$ be a regular or manageable
multiplicative unitary and $S= S_V$. Let $(\mu, \nu)$ be the $V_S$-covariant
representation defined by $V$. Then $\{({\rm id}\otimes\omega\circ\mu)(V_S):%
\omega\in {\cal L}(H)_*\}$ is dense in $\hat S_p$. On the other hand, if $S$
is an arbitrary Hopf $C^*$-algebra such that both $A_S$ and $A_{\hat S_p}$
are non-zero, the
density condition will be true for any $V_S$-covariant representation (see
the proof of [12, 5.3]).

\noindent (b) If the invariant mean $m$ in condition (iv) satisfies the
following stronger condition: $m(\{(\omega\circ\nu\otimes j)(V_S): \omega\in 
{\cal L} (H)_*\})\neq (0)$, then $S$ will have a Haar state (which is the
restriction of $m$) and is a compact quantum group.

\noindent (c) If the algebra $B$ in the above proposition has a bounded
right approximate identity, then $B$ has a bounded 2-sided approximate
identity (by Proposition 1.12(b)).

\medskip

\noindent {\bf Definition and Proposition 3.8:} Let $(T,X,S)$ be a general
Fourier duality. Then a $X$-covariant representation $(\mu,\nu)$ on $H$ is
said to be {\it amenable} if it satisfies one of the following equivalent
conditions\vspace{-3mm}:

\begin{enumerate}
\item[1.]  $T_0=\nu(T)$ is counital\vspace{-2mm}.

\item[2.]  There exists a net of unit vectors $\{\eta_i\}$ of $H$ such that
for any $\xi\in H$, $\|(\nu\otimes\mu)(X)(\eta_i\otimes\xi) -
\eta_i\otimes\xi\|$ converges to zero\vspace{-2mm}.

\item[3.]  $\nu^*({\cal L}(H)_*)$ has a bounded left approximate identity%
\vspace{-2mm}.
\end{enumerate}

\noindent {\bf Proof:} The proofs of (2) implies (3) and (3) implies (1) are
the same as those of Proposition 3.5 (note that $\nu^*({\cal L}(H)_*)$ is a
left Hopf ideal of $T^*$ and $\nu^*(T_0^*) = \overline{\nu^*({\cal L}(H)_*)}%
^{\sigma(T^*,T)}$). Suppose that (1) holds. $Y=(\nu\otimes {\rm id})(X)$ is
clearly a $T_0$-$S$-birepresentation such that $\{ ({\rm id}\otimes f)(Y):
f\in S^\#\}$ is dense in $T_0$. Therefore, by the proof of Lemma 1.6, the
counit of $T_0$ is a Hopf $*$-homomorphism $\pi$ of $T_0$ such that $%
(\pi\otimes {\rm id})(Y)=1$. It is clear that $\pi\circ\nu$ is weakly
contained in $\nu$ and condition (2) follows from the same argument as in
[5, 2.7.1].

\medskip

Since $S$ is dualizable, $T$ is a quotient Hopf $C^*$-algebra of $\hat S_p$
(by [12, 3.5]) and any $X$-covariant representation induces a $V_S$%
-covariant representation. In this case, the induced representation is
amenable if and only if the original one is.

\medskip

Moreover, we also have the corresponding results of Proposition 3.6 in the
case of a general Fourier duality. Note that $\{({\rm id}\otimes\omega\circ%
\mu)(X): \omega\in{\cal L}(H)_*\}$ is automatically dense in $T$ in the case
of Kac-Fourier duality or duality arising from a regular or manageable
multiplicative unitary (see Remark 3.7(a)).

\medskip

Next, we would like to study crossed products of coactions of amenable
Hopf $C^*$-algebras. We refer the readers to [12, \S 2 \& \S4] for the basic
definitions and properties of (full and reduced) crossed products. From now
on, $A$ and $B$ are two $C^*$-algebras and $\epsilon$ is a coaction of $S$
on $A$. Let $\epsilon^{\prime}$ be the coaction on $A\otimes B$ defined by $%
\epsilon^{\prime}(a\otimes b) = (\epsilon(a)\otimes b)^{\sigma_{23}}$ (where 
$\sigma_{23}$ is the flip of the second and the third variables). Using the
same sort of argument as in [11, 2.14 \& 3.2], we have the following
proposition.

\medskip

\noindent {\bf Proposition 3.9:} (a) If $(\mu,\nu)$ is an amenable $V_S$%
-covariant representation on $H$, then $A\times^{\mu,\nu}_\epsilon \hat S =
A\times_\epsilon \hat S$.

\noindent (b) For any $V_S$-covariant representation $(\mu,\nu)$, $%
(A\times_{\epsilon}^{\mu,\nu}\hat S)\otimes B = (A\otimes
B)\times_{\epsilon^{\prime}}^{\mu,\nu}\hat S$.

\noindent (c) If $A$ is nuclear, then $(A\times_\epsilon\hat
S)\otimes_{\max}B = (A\otimes B)\times_{\epsilon^{\prime}}\hat S$.

\medskip

This clearly implies the following result.

\medskip

\noindent {\bf Proposition 3.10:} Let $A$ be a nuclear Hopf $C^*$-algebra
with a coaction $\epsilon$ by an amenable Hopf $C^*$-algebra $S$. Then the
crossed product $A\times_\epsilon \hat S$ is nuclear. Consequently, if $S$
is a amenable Hopf $C^*$-algebra, then $\hat S_p$ is a nuclear $C^*$%
-algebra.

\medskip

The final statement applies, in particular, to the case when there exists a $%
V_S$-covariant representation and $\hat S_p$ is $H_2$-coamenable (see
Corollary 3.4).

\medskip

\noindent {\bf Remark 3.11:} We can also show that the corresponding results
of [11, 3.9-3.11] hold in this situation (with virtually the same
arguments). Hence we have the following: If $A$ is a $C^*$-exact $C^*$%
-algebra with coaction $\epsilon$ by an amenable Hopf $C^*$-algebra $S$,
then $A\times_\epsilon \hat S$ is $C^*$-exact.

\medskip

Finally, we would like to see the relation between the $H_i$-coamenability ($%
i=1,2,3$) and amenability in the case of Kac-Fourier duality.

\medskip

\noindent {\bf Proposition 3.12:} Let $(T,X,S)$ be a Kac-Fourier duality.
Then the following statements are equivalent.

\noindent (i) $T$ is left (or right) $H_1$-coamenable.

\noindent (ii) $T$ is counital.

\noindent (iii) $S$ is amenable.

\noindent In this case, $T\cong\hat S_p$. Consequently, if $S$ is
amenable, then $A_S$ has a right invariant mean.

\noindent {\bf Proof:} Let $(\mu,\nu)$ be a $X$-covariant representation
such that both $\mu$ and $\nu$ are injective (see Definition 1.5(b)). By
Proposition 3.8 and Theorem 1.8, (i) is equivalent to (ii). It is clear from
Proposition 3.5 that (ii) implies (iii). Moreover, for any $V_S$-covariant
representation $(\tilde\mu,\tilde\nu)$ (note that $S$ is dualizable by
assumption and $(\hat S_p,V_S,S)$ is the intrinsic duality of $S$), $%
\tilde\nu(\hat S_p) = T$ (by [12, 5.13]) and we have (iii) implies (ii) (by
Proposition 3.5 again). The last two statements are clear (using Proposition
3.6).

\medskip

\noindent {\bf Corollary 3.13:} Let $(T,X,S)$ be a Kac-Fourier duality.

\noindent (a) The one sided and two sided $H_i$-coamenability (for $i=1,2,3$)
of $T$ are all the same.

\noindent (b) If $T$ is $H_i$-coamenable for some $i=1,2,3$, then $T$ is a
nuclear $C^*$-algebra.

\noindent {\bf Proof:} Note that by Theorem 1.8, $A^l_T = A_T = A^r_T$.
Moreover, Proposition 3.12 and Proposition 1.12(b) show that the existence
of bounded left (or right or 2-sided) approximate identity on $A_S$ is
equivalent to $T$ being counital. Hence, part (a) is proved. Now part (b)
clearly follows from part (a), Propositions 3.12 and 3.10.

\medskip

This initiates us to make the following definition.

\medskip

\noindent {\bf Definition 3.14:} A Kac-Fourier duality $(T,X,S)$ is said to
be {\it amenable} (respectively, {\it coamenable}) if $S$ (respectively, $T$%
) is amenable.

\medskip

\noindent {\bf Remark 3.15:} (a) For a Kac-Fourier duality $(T,X,S)$, 
$H_1$-coamenability and $H_2$-coamenability of $\hat S_p$ are equivalent to each
other. In fact, we have $A_{\hat S_p} = A_T$ and if $\hat S_p$ is 
$H_2$-coamenable, then any $V_S$-covariant representation is amenable (by Corollary
3.4) and so by Proposition 3.5, $A_{\hat S_p}$ have a bounded approximate
identity.

\noindent (b) However, we do not know whether the $H_3$-coamenability of $\hat
S_p$ is equivalent to the other two. Note that it is the case for locally
compact groups (see [9, Theorem 1]).

\medskip

\medskip

\medskip

\noindent {\bf REFERENCE:}

\medskip

\medskip

\noindent [1] S. Baaj, Repr\'{e}sentation r\'{e}guli\`{e}re du groupe
quantique des d\'{e}placements de Woronowicz, Asterisque (1995), no. 232,
11-48.

\noindent [2] S. Baaj and G. Skandalis, Unitaires multiplicatifs et
dualit\'e pour les produits crois\'es de $C^{*}$-alg\`ebres, Ann. scient.
\'Ec. Norm. Sup., $4^{e}$ s\'erie, t. 26 (1993), 425-488.

\noindent [3] \'{E}. Blanchard, D\'{e}formations de $C\sp *$-alg\`{e}bres de
Hopf, Bull.-Soc.-Math.-France 124 (1996), no. 1, 141--215.

\noindent [4] Jean De Canni\`{e}re and Ronny Rousseau, The Fourier Algebra
as an Order Ideal of the Fourier-Stieltjes Algebra, Math. Z. 186 (1984),
501-507.

\noindent [5] M. Enock and J.-M. Schwartz, Algebres de Kac moyennables,
Pacific J. Math. 125 (1986), 363-379.

\noindent [6] J.M.G. Fell, Weak containment and induced representations of
groups, Cand. J. Math. 14 (1962), 237-268.

\noindent [7] B. Forrest, Amenability and Derivations of the Fourier
algebra, Proc. Amer. Math. Soc. 104 (1988), 437-442.

\noindent [8] A.T. Lau, Fourier and Fourier-Stieltjes algebras of a locally
compact group and amenability (Topological vector spaces, algebras and
related areas), Pit. Res. Notes in Math. 316 (1994), 79-92.

\noindent [9] V. Losert, Properties of the Fourier algebras that are
equivalent to amenability, Proc. Amer. Math. Soc. 92 (1984), 347-354.

\noindent [10] C. Nebbia, Multipliers and asymptotic behaviour of the
Fourier algebra of nonamenable groups, Proc. Amer. Math. Soc. 84 (1982),
549-554.

\noindent [11] C.K. Ng, Coactions and crossed products of Hopf $C^*$%
-algebras, Proc. Lond. Math. Soc.(3) 72 (1996), 638-656.

\noindent [12] C.K. Ng, Duality of Hopf $C^*$-algebras, preprint.

\noindent [13] T.W. Palmer, Banach algebras and the general theory of $*$%
-algebras Vol. I: Algebras and Banach algebras, Encyclopedia of Math. (Camb.
Univ. Press, 1994).

\noindent [14] A.T. Paterson, Amenability, Math. Surveys and Monographs 29
(Amer. Math. Soc., 1988).

\noindent [15] J.-P. Pier, Amenable Banach algebras, Pit. Res. Notes in
Math. 172 (Longman, 1988).

\noindent [16] Z.J. Ruan, The operator amenability of $A(G)$, Amer. J. of
Math. 117 (1995), 1449-1474.

\noindent [17] J.C.S. Wong and A. Riazi, Characterisations of amenable
locally compact semigroups, Pacific J. of Math. 108 (1983), 479-496.

\noindent [18] S. L. Woronowicz, From multiplicative unitaries to quantum
groups, Internat. J. Math. 7(1996), no.1, 127-149.

\noindent [19] S. L. Woronowicz, Compact quantum groups, 
Sym\'{e}tries quantiques (Les Houches, 1995), North-Halland, Amsterdam 1998, 
pp. 845-884.

\medskip

\medskip

\medskip

\noindent Mathematical Institute, 24-29 St. Giles, Oxford OX1 3LB, United
Kingdom.

\medskip \noindent $e$-mail address: ng@maths.ox.ac.uk

\end{document}